\documentclass[12pt]{amsart}

\usepackage{etex}
\usepackage{amscd}
\usepackage{amsmath}
\usepackage{amsfonts}
\usepackage{amssymb}
\usepackage{graphics}
\usepackage{epsfig}
\usepackage{pictex} % this is only necessary for the \setshagrid and related commands
\usepackage{pict2e} % для рисования окружностей большого радиуса, без этого пакета радиус ограничен
\usepackage{url}

\font\Bbb=msbm10 scaled \magstep 2

\def\C{\hbox{\Bbb C}}

\newtheorem{theorem}{\bf Theorem}
\newtheorem{lemma}[theorem]{\bf Lemma}

\newtheorem{example}[theorem]{\bf Example}

\newtheorem{corollary}[theorem]{\bf Corollary}
\newtheorem{definition}[theorem]{\bf Definition}

\title[A Graceful Basis]
{A Graceful Basis in the Solution Space of an ODE with Constant Coefficients}

\author{Timur Sadykov}
\address{Plekhanov Russian University of Economics
\newline \indent 115054, Moscow, Russia}
\email{Sadykov.TM@rea.ru}
\thanks{This research was performed in the framework of the state task in the field of scientific activity
of the Ministry of Science and Higher Education of the Russian Federation, project no. FSSW-2023-0004.}

\begin{document}

\begin{abstract}
We~revisit the classical problem of construction of a fundamental system of solutions to a linear ODE
whose elements remain analytic and linearly independent for all values of the roots of the characteristic polynomial.
\end{abstract}

\medskip

\maketitle

It~is well-known that the canonical basis in the space of analytic solutions to the ordinary differential equation with constant coefficients
\begin{equation}
\label{ODEWithConstCoeff}
\left(\prod_{j=1}^{m} \left( \frac{d}{dx} - \alpha_j \right)\right) f(x) = 0
\end{equation}
has the form
\begin{equation}
\label{ODEWithConstCoeffStandardBasis}
\left\{
e^{\alpha_j x}, \, j=1,\ldots,m
\right\}
\end{equation}
in the case when all the roots~$\alpha_j$ of its characteristic polynomial are different.
For multiple roots, the exponents~(\ref{ODEWithConstCoeffStandardBasis}) become linearly dependent,
but the space of analytic solutions to~(\ref{ODEWithConstCoeff}) has dimension~$m$
no matter what the values of $\alpha:=(\alpha_1,\ldots,\alpha_m)\in\C^m$ are.

The search for fundamental systems of solutions to~(\ref{ODEWithConstCoeffStandardBasis}) that do not degenerate for confluent values of
its parameters $\alpha\in\C^m$ (i.e., such that at least two of them are equal) received a lot of attention,
both classically and recently (see~\cite[\S~26]{Arnold-ODE} and the references therein).
The canonical method of multiplying an exponent in~(\ref{ODEWithConstCoeffStandardBasis}) with an arbitrary polynomial
whose degree is one smaller than the multiplicity of a root yields a basis whose elements depend on~$\alpha$ in a discontinuous fashion.
Besides, it does not work in the case of symbolic~$\alpha_j$ whose values are not a priori known.
In~\cite{Gatto-Shcherbak}, a universal fundamental system of solutions to a generic ordinary linear differential equation with constant coefficients
is given in terms of formal power series.
In~\cite[Section~10.1]{Sturmfels}, the space of solutions to a second-order differential equation is endowed with a basis
whose elements remain analytic and linearly independent for all values of the roots of the characteristic polynomial.
This behavior of a basis is referred to in~\cite[p.~134]{Sturmfels} as "graceful".
Following~\cite{Sturmfels}, we adopt the next definition.

\begin{definition}
\label{def:gracefulBasis}\rm
An $m$-tuple $f_j(x;\alpha),\, j=1,\ldots,m$ of analytic solutions to the equation~(\ref{ODEWithConstCoeffStandardBasis})
is called {\it a graceful basis} if the following hold:

1) $f_j(x;\alpha)$ are entire functions of $x\in\C$ and $\alpha\in\C^m;$

2) for any fixed $\alpha\in\C^m$ the univariate functions $f_j(x;\alpha)$ are linearly independent.
\end{definition}

Clearly the canonical basis~(\ref{ODEWithConstCoeffStandardBasis}) is not graceful, nor is it its standard completion by polynomial factors
in the case of confluent values of the characteristic roots~$\alpha_1,\ldots,\alpha_m.$ The next lemma remedies this flaw.

\begin{lemma}
\label{lemma:gracefulBasis}
Denote by $\mathcal{S}_d (\alpha)$ the elementary symmetric polynomial of degree~$d$ with the variables~$\alpha$
and by~$[j]$ the omission of the index~$j$ in a list.
The family of~$m$ functions
\begin{equation}
\label{gracefulBasis}
\mathfrak{g}_{i}(x):=
\sum_{j=1}^m
\frac{(-1)^{m-i} \, \mathcal{S}_{m-i}(\alpha_1,\ldots [j] \ldots,\alpha_m)}
{(\alpha_j-\alpha_1)\cdot\ldots [j] \ldots\cdot(\alpha_j-\alpha_m)
}
\, e^{\alpha_j x}, \quad i=1,\ldots,m
\end{equation}
is a graceful basis in the space of analytic solutions to~(\ref{ODEWithConstCoeff}).
\end{lemma}
For $m>1$ the basis~(\ref{gracefulBasis}) is different from the universal basis given in~\cite[formula~(4)]{Gatto-Shcherbak},
although the sum of its initial series equals~$\mathfrak{g}_{m}(x)$.
For the maximally confluent values of the roots of the characteristic polynomial, i.e.,
for $\alpha_1 = \ldots = \alpha_m = 0,$ the basis~(\ref{gracefulBasis}) turns into the family of monomials
$1,x,\frac{x^2}{2},\ldots,\frac{x^{m-1}}{(m-1)!}.$
\begin{proof}

Straightforward expansion of~(\ref{gracefulBasis}) into power series shows that the hyperplanes $\alpha_j=\alpha_k$
are removable singularities of~$\mathfrak{g}_{i}(x)$ which are thereby entire functions of $x\in\C$ and $\alpha\in\C^m.$
Alternatively, this can be seen by means of the residue theory.
For instance, the last element~(\ref{gracefulBasis}) of admits the representation
$$
\mathfrak{g}_{m}(x) =
-\underset{t=\infty}{\rm res}\, \frac{e^{t x}}{(t-\alpha_1)\ldots(t-\alpha_m)} =
\frac{1}{2\pi\sqrt{-1}}\int\limits_{|t|=R} \frac{e^{t x} \, dt}{(t-\alpha_1)\ldots(t-\alpha_m)},
$$
with $R>\max\limits_{j} |\alpha_j|.$ The above integral is clearly an entire function of~$x$ and~$\alpha.$

Denote by $V(\alpha)$ the Vandermonde matrix defined by the roots of the characteristic polynomial of~(\ref{ODEWithConstCoeff}):
$$
V(\alpha)\equiv
V(\alpha_1,\ldots,\alpha_m):=
\left(
\begin{array}{ccccc}
1      & \alpha_{1} & \ldots & \alpha_{1}^{m-1} \\
1      & \alpha_{2} & \ldots & \alpha_{2}^{m-1} \\
\ldots & \ldots     & \ldots &   \ldots         \\
1      & \alpha_{m} & \ldots & \alpha_{m}^{m-1} \\
\end{array}
\right).
$$
It~is well-known that the determinant of this matrix equals $\prod\limits_{1\leq j<k \leq m} (\alpha_j - \alpha_k)$.
Since~(\ref{gracefulBasis}) are linear combinations of the elements of the canonical basis~(\ref{ODEWithConstCoeffStandardBasis}),
they satisfy the equation~(\ref{ODEWithConstCoeff}).
By~\cite[p.~38]{Knuth} (see also~\cite{Macon-Spitzbart,Tou} and the references therein),
the coefficients of these linear combinations form the inverse to the matrix~$V(\alpha)$ defined by nonconfluent characteristic roots.
By the conservation principle for analytic differential equations~\cite{Marichev} the functions~(\ref{gracefulBasis}) solve~(\ref{ODEWithConstCoeff})
for all~$\alpha\in\C^m.$

It~remains to show that the functions~(\ref{gracefulBasis}) are linearly independent.
Denote $|\alpha|:=\alpha_1 + \ldots + \alpha_m.$
The Wronskian determinant of the basis~(\ref{ODEWithConstCoeffStandardBasis}) equals
$e^{|\alpha| x} \prod\limits_{1\leq j<k \leq m} (\alpha_j - \alpha_k)$
and vanishes if and only if the characteristic polynomial of~(\ref{ODEWithConstCoeff}) has multiple roots.
It~therefore follows that for nonconfluent~$\alpha$ the Wronskian determinant of the family of functions~(\ref{gracefulBasis}) is given~by
$$
{\rm Wr}\left(\mathfrak{g}_{1}(x),\ldots,\mathfrak{g}_{m}(x); \, x\right) =
\det{(V(\alpha))^{-1}} \, {\rm Wr} \left(e^{\alpha_1 x},\ldots,e^{\alpha_m x}; \, x\right) =
e^{|\alpha|x}.
$$
By~the uniqueness theorem for analytic functions this equality holds for all $\alpha\in\C^m.$
Thus the above Wronskian does not vanish anywhere and hence~(\ref{gracefulBasis})
is indeed a graceful basis in the space of solutions to~(\ref{ODEWithConstCoeff}).
\end{proof}

Skipping the explicit inverse to the nondegenerate Vandermonde matrix, we arrive at the following particularly simple form of~(\ref{gracefulBasis}):

\begin{corollary}
\label{cor:simpleFormulaForGracefulBasis}
The functions $\left(V(\alpha)\right)^{-1} \left(e^{\alpha_1 x},\ldots, e^{\alpha_m x}\right)^T$
form a graceful basis in the space of analytic solutions to~(\ref{ODEWithConstCoeff}).
\end{corollary}

\begin{example}
\label{ex:m=2,3}\rm
A~graceful basis in the space of solutions to the linear homogeneous ordinary differential equation of order two
$\left(\frac{d}{dx} - \alpha_1 \right)\left( \frac{d}{dx} - \alpha_2 \right)f(x)=0$
has the following form (cf.~\cite[Section~10.1,~(10.5)]{Sturmfels}):
$$
\left\{\frac{-\alpha_2 e^{\alpha_1 x} + \alpha_1 e^{\alpha_2 x}}{\alpha_1-\alpha_2},\right. \quad
\left.\frac{e^{\alpha_1 x} - e^{\alpha_2 x}}{\alpha_1-\alpha_2}\right\}.
$$

The entire functions
$$
\left\{
\frac{\alpha_2 \alpha_3 e^{\alpha_1 x}}{\left(\alpha_1-\alpha_2\right) \left(\alpha_1-\alpha_3\right)}-
\frac{\alpha_1 \alpha_3 e^{\alpha_2 x}}{\left(\alpha_1-\alpha_2\right) \left(\alpha_2-\alpha_3\right)}+
\frac{\alpha_1 \alpha_2 e^{\alpha_3 x}}{\left(\alpha_1-\alpha_3\right) \left(\alpha_2-\alpha_3\right)},
\right.
$$
$$
-\frac{\left(\alpha_2+\alpha_3\right) e^{\alpha_1 x}}{\left(\alpha_1-\alpha_2\right) \left(\alpha_1-\alpha_3\right)}
+\frac{\left(\alpha_1+\alpha_3\right) e^{\alpha_2 x}}{\left(\alpha_1-\alpha_2\right) \left(\alpha_2-\alpha_3\right)}
-\frac{\left(\alpha_1+\alpha_2\right) e^{\alpha_3 x}}{\left(\alpha_1-\alpha_3\right) \left(\alpha_2-\alpha_3\right)},
$$
$$
\left.
\frac{e^{\alpha_1 x}}{\left(\alpha_1-\alpha_2\right) \left(\alpha_1-\alpha_3\right)}-
\frac{e^{\alpha_2 x}}{\left(\alpha_1-\alpha_2\right) \left(\alpha_2-\alpha_3\right)}+
\frac{e^{\alpha_3 x}}{\left(\alpha_1-\alpha_3\right) \left(\alpha_2-\alpha_3\right)}
\right\}
$$
form a graceful basis in the space of solutions to~(\ref{ODEWithConstCoeff}) for $m=3.$
\end{example}

%------------------------------------------------------------------------------------------------------------

\end{document}